\documentclass[letterpaper, 10 pt, conference]{IEEEtran}
\usepackage[utf8]{inputenc}
\usepackage[english]{babel}
\usepackage[backend=bibtex]{biblatex}
\bibliography{bibliography.bib}

\usepackage{lettrine}
\usepackage{amsmath}
\usepackage{booktabs}
\usepackage{amssymb,epsfig,multirow}
\usepackage{oldlfont}
\usepackage{algorithm}

\usepackage{algpseudocode}
\makeatletter
\def\BState{\State\hskip-\ALG@thistlm}
\makeatother
\usepackage{graphicx}
\usepackage{url}
\usepackage{comment}
\usepackage{float}
\usepackage{color}
\usepackage{enumerate}
\usepackage{version} 
\usepackage{algorithmicx}
\usepackage{tikz}
\usetikzlibrary{arrows.meta}
\usepackage{tkz-euclide}
\usepackage{siunitx}
\usepackage{subcaption}
\usepackage{algpseudocode}
\PassOptionsToPackage{noend}{algpseudocode}
\algtext*{EndWhile}
\algtext*{EndIf}
\algtext*{EndFor}
\errorcontextlines\maxdimen

\makeatletter
\newcommand*{\algrule}[1][\algorithmicindent]{\makebox[#1][l]{\hspace*{.5em}\vrule height .75\baselineskip depth .25\baselineskip}}%

\newcount\ALG@printindent@tempcnta
\def\ALG@printindent{%
    \ifnum \theALG@nested>0
        \ifx\ALG@text\ALG@x@notext
            \addvspace{-3pt}
        \else
            \unskip
            \ALG@printindent@tempcnta=1
            \loop
                \algrule[\csname ALG@ind@\the\ALG@printindent@tempcnta\endcsname]%
                \advance \ALG@printindent@tempcnta 1
            \ifnum \ALG@printindent@tempcnta<\numexpr\theALG@nested+1\relax
            \repeat
        \fi
    \fi
    }%
\usepackage{etoolbox}
\patchcmd{\ALG@doentity}{\noindent\hskip\ALG@tlm}{\ALG@printindent}{}{\errmessage{failed to patch}}
\makeatother

\usepackage{pifont}
\algdef{SE}[DOWHILE]{Do}{doWhile}{\algorithmicdo}[1]{\algorithmicwhile\ #1}%
\usepackage{verbatim}

\usepackage{times}
\newtheorem{theorem}{Theorem}
\numberwithin{theorem}{section}

\newtheorem{prop}[theorem]{Proposition}

\usetikzlibrary{shapes,snakes}

\newtheorem{defn}{Definition}[section]

\newcommand{\C}{{\mathcal{C}}}
\newcommand{\E}{{\mathcal{E}}}
\newcommand{\A}{{\mathcal{A}}}
\newcommand{\Ll}{{\mathcal{L}}}
\newcommand{\N}{{\mathcal{N}}}
\newcommand{\B}{{\mathcal{B}}}

\newcommand{\Ss}{{\mathcal{S}}}

\title{ \bigskip
\medskip
\LARGE \bf Local Optimization of MAPF solutions on Directed Graphs}

\vspace{-3pt}

\author{S. Ardizzoni$^{1}$, I. Saccani$^{1}$, L. Consolini$^{1}$,  M. Locatelli$^{1}$  \\
       { \small $^{1}$Dipartimento di Ingegneria e Architettura, Universit\`a di Parma, Parco Area delle Scienze, 181/A, Parma, Italy}

\vspace{-3pt}       
       
}

\date{}
\begin{document}
	\tikzstyle{abstract}=[rectangle, draw=black, rounded corners, fill=blue!40, drop shadow,
	text centered, anchor=north, text=white, text width=3cm]
	\tikzstyle{comment}=[rectangle, draw=black, rounded corners, fill=red, drop shadow,
	text centered, anchor=north, text=white, text width=3cm]
	\tikzstyle{myarrow}=[->, >=open triangle 90, thick]
	\tikzstyle{line}=[-, thick]
	
	\maketitle
	
	\begin{abstract}                          
Among sub-optimal MAPF solvers, rule-based algorithms are particularly appealing since they are complete. Even in crowded scenarios, they allow finding a feasible solution that brings each agent to its target, preventing deadlock situations. However, generally, rule-based algorithms provide solutions that are much longer than the optimal one. 
The main contribution of this paper is the introduction of an iterative local search procedure in MAPF. We start from a feasible suboptimal solution and we perform a local search in a neighborhood of this solution, to find a shorter one. Iteratively, we repeat this procedure until the solution cannot be shortened any longer.
At the end, we obtain a solution, that is still sub-optimal, but, in general, of much better quality than the initial one.
We use dynamic programming for the local search procedure. Under this respect, the fact that our search is \emph{local} is fundamental to reduce the time complexity of the algorithm. Indeed, if we apply a standard dynamic programming the number of explored states grows exponentially with the number of agents. As we will see, the introduction of a locality constraint allows solving the (local) dynamic programming problem in a time that grows only polynomially with respect to the number of agents.
	\end{abstract}

\section{Introduction}
We focus on the Multi-Agent Path Finding (MAPF) problem. 
We consider a directed graph and a set of agents. Each agent occupies a different node and can move to free nodes, i.e., nodes not occupied by other agents.
The MAPF problem consists in computing a sequence of movements that repositions all agents to assigned target nodes, avoiding collisions. 
The main motivation comes from the management of fleets of automated guided vehicles (AGVs). AGVs move items between different locations in a warehouse. Each AGV follows predefined paths, that connect the locations in which items are stored or processed. We associate the layout of the paths to a directed graph.
The nodes represent positions in which items are picked up and delivered,
together with additional locations used for routing. The directed arcs represent the precomputed paths that connect these locations. If various AGVs move in a small scenario, each AGV represents an obstacle for the other ones. In some cases, the fleet can reach a deadlock situation, in which every vehicle is unable to reach its target. Hence, it is important to find a feasible solution to MAPF, even in crowded configurations.

{\bf Literature review.}
Various works address the problem of finding the optimal solution of MAPF (i.e., the solution with the minimum number of moves). For instance, Conflict Based Search (CBS) is a two-level algorithm which uses a search tree, based on 
conflicts between individual agents (see~\cite{sharon2015}).
However, finding the optimal solution of MAPF is NP-hard (see~\cite{yu2013}), and computational time grows exponentially with the number of agents. Therefore, optimal solvers are usually applied when the number of agents is relatively small.
Conversely, sub-optimal solvers are usually employed when the number of agents is large. In such cases, the aim is to quickly find a path for the different agents, and it is often intractable to guarantee that a given solution is optimal. Among these, search-based solvers aim to provide a high quality solution, but are not complete (i.e., they are not always able to return a feasible solution). A prominent example is Hierarchical Cooperative A$^*$ (HCA$^*$) \cite{silver2021}, in which agents are planned one at a time according to some predefined order.
Instead, rule-based approaches include specific movement rules for different scenarios. They favor completeness at low computational cost over solution quality. One of the first important results in this field is from Kornhauser's thesis \cite{kornhauser1984}, which presents a rule-based procedure to solve MAPF (or to establish that MAPF has no feasible solution). Kornhauser shows that solutions found with his method have cubic length complexity. However, to our knowledge, this procedure has never been fully implemented. Two relevant recent rule-based algorithms are TASS \cite{khorshid2021} and \textit{Push and Rotate} \cite{wilde2014} \cite{alotaibi2016}. TASS is a tree-based agent swapping strategy which is complete on every tree, while \textit{Push and Rotate} solves every MAPF instance on graphs that contains at least two holes (i.e., unoccupied vertices). 
Reference~\cite{krontiris2021} presents a method that converts the graph into a tree (as in~\cite{goraly2010}), and solves the resulting problem with TASS. 
Rule-based algorithms are also used for directed graph with at least two unoccupied nodes. In particular, reference \cite{ardizzoni2022} presents diSC algorithm, which solve any MAPF instance on strongly connected digraphs, i.e., directed graphs
in which it is possible to reach any node starting from any
other node.
Another relevant reference is~\cite{botea2018}, which solves MAPF on the specific class of biconnected digraphs, i.e., strongly connected digraphs where the undirected graphs obtained by ignoring the edge orientations have no cutting vertices. The proposed diBOX algorithm has polynomial complexity with respect to the number of nodes.

{\bf Motivations and statement of contribution.}

Among sub-optimal MAPF solvers, rule-based algorithms are particularly appealing since they are complete. Even in crowded scenarios, they allow finding a feasible solution that brings each agent to its target, preventing deadlock situations. However, generally, rule-based algorithms provide solutions that are much longer than the optimal one. This is a crucial limitation in industrial applications.




The main contribution of this paper is the introduction of an iterative local search procedure in MAPF. We start from a feasible suboptimal solution, for instance the one provided by a rule-based algorithm. We perform a local search in a neighborhood of this solution, to find a shorter one. Iteratively, we repeat this procedure until the solution cannot be shortened any longer.
At the end, we obtain a solution, that is still sub-optimal, but, in general, of much better quality than the initial one.

We use dynamic programming for the local search procedure. Under this respect, the fact that our search is \emph{local} is fundamental to reduce the time complexity of the algorithm. Indeed, in principle, it is possible to solve the general MAPF problem by dynamic programming. However, the number of explored states grows exponentially with the number of agents, so that we cannot apply standard dynamic programming to problems involving many agents.  As we will see, the introduction of a locality constraint allows solving the (local) dynamic programming problem in a time that grows only polynomially with respect to the number of agents (see Theorem~\ref{thm_complx}).

\section{Problem Definition}

\subsection{MAPF problems}

Let $G=(V,E)$ be a directed graph, with vertex set $V$ and edge set $E$.
We assign a unique label to each agent, and set $P$ contains these labels. 
A \textit{configuration} is a function  $\A:P \rightarrow V$ that assigns the occupied vertex to each agent. A configuration is \textit{valid} if it is injective (i.e., each vertex is occupied by at most one agent). Set $\C \subset \{P \to V\}$ represents all valid configurations. 




		
		
		
	
	Time is assumed to be discretized. At every time step, each agent occupies one vertex and executes a single action. 
	There are two types of actions: \emph{wait} and \emph{move}. We denote the wait action by $\iota$. An agent that executes this action remains in its current vertex for another time step. We denote a move action by $u \rightarrow v$.  In this case, the agent moves from its current vertex $u$ to an adjacent vertex $v$  (i.e., $(u,v) \in E$). Therefore, the set of all possible actions for a single agent is
	
	\[ \bar{E} = E \cup \{\iota\}.\]

	Function $\rho: \C \times \bar{E} \rightarrow \C$ is a partially defined transition function such that $\A' =\rho(\A,u \rightarrow v)$  is the configuration obtained by moving an agent from $u$ to $v$:
	
	\begin{equation} 
		\label{c}
		\A'(q):=  \Bigg\{
		\begin{array}{ll}       	
			
			v, & \text{if } \A(q)  = u ;\\

			\A(q), & \text{otherwise } .\\
		\end{array}
	\end{equation}
	   Notation $\rho(\A,u \rightarrow v)!$ means that the function is well-defined. In other words
	$\rho(\A, u \rightarrow v)!$ if and only if $(u,v)\in E$ and $\A' \in \C$.  Moreover, $(\forall \A \in \C)\; \rho(\A,\iota)!$ and $\rho(\A,\iota)=\A$.
	
	Since the movements of the agents can be synchronous, at each time step an action is an element of $\E=\bar{E}^{|P|}$, $a=(a_1,\dots ,a_{|P|})$ where $a_i$ is the single move of agent $i$.
	We can extend function $\rho:\C \times \bar{E} \to \C$ to $\rho: \C \times \E  \to \C$, by setting  $\A'=\rho(\A,  a) $ equal to the configuration obtained by moving agent $i$ along edge $a_i$ (or by not moving the agent if $a_i=\iota$).
	In this case,
	 $(\forall a \in \E, \A \in \C)$  $\rho(\A,  a) !$ if and only if the following conditions hold:
	 
	 \begin{enumerate}
	 	\item $\A' \in \C$: two ore more agents cannot occupy the same vertex at the same time step;
 	     \item $\forall i=1,\dots, |P|$, if $a_i=(u,v)$, then $\not \exists j \in \{1,\dots, |P|\}$ such that $a_j=(v,u)$: two agents cannot swap locations in a single time step.
	 \end{enumerate} 
	
	
	We represent plans as ordered sequences of actions.
	It is convenient to view the elements of $\E$ as the symbols of a language.
	We denote by $\E^*$ the Kleene star of $\E$, that is the set of
	ordered sequences of elements of $\E$ with arbitrary length, together with the empty string $\epsilon$:
	\[
	\E^*=\bigcup_{i=1}^\infty \E^i \cup \{\epsilon\}.
	\]
	We extend function $\rho:\C \times \E \to \C$ to $\rho: \C \times \E^*  \to \C$.  $(\forall s \in \E^*, e \in \E, \A \in \C)$  $\rho(\A,  se) !$ if and only if $\rho(\A, s)!$ and $\rho(\rho(\A, s), e)!$ and, if $\rho(\A, se)!$, then $\rho(\A, se)=\rho(\rho(\A, s),e)$.

	Note that $\epsilon$ is the trivial plan that keeps all agents and holes at their positions.\\

	We denote by $\E^*_{\A}$ the set of plans such that $\rho(\A,  f) $ is well defined:
	
	\[\E^*_{\A} = \{ f \in \E^* : \, \rho(\A,  f)!\}. \]

	

	The problem of detecting a feasible solution is the following:

	\begin{defn}{(\textbf{Feasibility MAPF problem}).}
		Given a digraph $G = (V,E)$, an agent set $P$, an initial valid configuration $\A^s$, and a final valid configuration $\A^t$, find a plan $f$ such that $\A^t = \rho(\A^s,f)$.\\
	\end{defn}

Now, for a feasible plan $f$, we define $|f|$ as the length of plan $f$, i.e., the number of time steps needed to let all agents reach the final configuration through plan $f$. Furthermore, given $k \in  {\mathbb{N}} $, we denote by $f_k$ the $k$-th prefix of $f$ (that is, the prefix of $f$ of length $k$, made up of the first $k$ actions of $f$). Note that $|f_k|=k$.

We aim to solve a given MAPF instance while minimizing a global cumulative cost function. We employ the cost function called \textit{Makespan}, equal to the time when the last agent reaches its destination (i.e., the maximum of the individual costs). 

\begin{defn}{(\textbf{Optimization MAPF problem}).}
Given $\A^s$ and $\A^t$ initial and final valid configurations on a digraph $G$, the optimization MAPF problem with \textit{Makespan} is defined as
\begin{equation} 
	\label{op}
	\centering 
	\begin{array}{ll}       	
		
	\min{|f|} & \\
	
	      & \\
		
	\quad 	\text{s.t.} \;\; \A^t = \rho(\A^s,f) &  \\
	
	  & \\
	
	\quad   \quad  \quad f \in \E^*_{\A^s}. & \\
	
	\end{array}
\end{equation}
\end{defn}

Other cost functions have also been used in the literature. \textit{Sum-of-costs}, for example, is the summation, over all agents, of the number of time steps that an agent employs to reach its target without leaving it again. Unfortunately, finding the optimal solution, i.e., the minimal \textit{Makespan} or \textit{sum-of-costs}, has been shown to be NP-hard \cite{yu2013}. Therefore, in this paper we propose an approach to detect a good quality sub-optimal solution in polynomial time.

\subsection{Distances}
\label{sec:distances}
As said, we propose a solution approach based on the exploration of a \emph{neighborhood} of a reference plan. To define a neighborhood, we introduce distances between vertices, configurations, and plans. 
Let $G=(V,E)$ be  a digraph and $P$ be a set of agents. We define the distance of vertex $u$ from vertex $v$ as the length of the shortest path on $G$ from $v$ to $u$:

\[ d: V \times V \rightarrow  {\mathbb{N}} \quad \quad d(u,v) = \ell(\pi_{vu}),\]

where $\pi_{vu}$ is the shortest path in $G$ from $v$ to $u$ and $\ell(\pi_{vu})$ is the length of that path, defined as the number of edges of $\pi_{vu}$. Note that $d$ is not symmetrical, since $\pi_{uv}$ and $\pi_{vu}$ can be different. Next, we define the distance of configuration $\A^1$ from configuration $\A^2$ as the sum of the distances between the vertices that each agent occupies in the two configurations:
\[ d: \C \times \C \rightarrow  {\mathbb{N}} \quad \quad d(\A^1,\A^2) = \sum_{p \in P} d(\A^1(p),\A^2(p)).\]

Finally, we define the asymmetrical distance between two plans in $\E^*$. To do that, we associate to each plan a function in $\Ll = \{ \psi: {\mathbb{N}} \rightarrow \C \}$ using the following 

\[\begin{array}{l} \Phi_{\A} : \E^*_{\A} \rightarrow \Ll \quad \quad  \\
	\\
	
\quad \quad \quad \quad 	f \rightarrow \psi_f(k):= \bigg\{ \begin{array}{ll}\rho(\A,f_k), & k < |f|, \\
\rho(\A,f), & k \geq |f|,
  \end{array}

 \end{array}\]

where $\psi_f$ is the function which associates to each $k \in {\mathbb{N}}$ the configuration at step $k$, that depends on the $k$-th prefix of $f$. We define the distance of plan $f$ from plan $g$ as the distance between the associated functions $\Phi_{\A}(f),\Phi_{\A}(g)$:

\[d: \E^*_{\A} \times \E^*_{\A} \rightarrow {\mathbb{N}} \quad \quad  d(f,g):= \bar{d}( \Phi_{\A}(f), \Phi_{\A}(g)).\]

We can define $\bar{d}$ in various ways, leading to different definitions of the distance between plans $f$ and $g$:

\begin{enumerate}
\item $\infty$-distance:
\[\bar{d}_{\infty}( \Phi_{\A}(f), \Phi_{\A}(g)) := \max_{1 \leq k \leq \min \{|f|,|g|\}}\; d(\psi_f(k), \psi_g(k)) ;\]
	
	\item $1$-distance: 
		\[ \bar{d}_{1}( \Phi_{\A}(f), \Phi_{\A}(g)) := \sum_{k =1}^{ \min \{|f|,|g|\}}\; d(\psi_f(k), \psi_g(k)) ;\]

\item max-min distance:
\[ \bar{d}_{\infty}^*( \Phi_{\A}(f), \Phi_{\A}(g)) :=\max_{k \in {\mathbb{N}}}\; \min_{h \in {\mathbb{N}}}d(\psi_f(k), \psi_g(h)) ;\]

\item sum-min distance:
\[\bar{d}_{1}^*( \Phi_{\A}(f), \Phi_{\A}(g)):= \sum_{k =1}^{ \min \{|f|,|g|\}} \min_{h \in {\mathbb{N}}} d(\psi_f(k), \psi_g(h)).\]
	
\end{enumerate}

Namely, with the $\infty$-distance, the distance of plans $f$ and $g$ corresponds to the maximum, with respect to time-step $k$, of the distance between the corresponding configurations at $k$. With the $1$-distance, this distance corresponds to the sum, with respect to time-step $k$, of the distances between the corresponding configurations at $k$. With the max-min distance (respectively, the sum-min distance), this distance corresponds to the maximum (respectively, the sum) with respect to $k$, of the distance of the configuration that plan $f$ reaches at step $k$ with respect to the set of all configurations encountered by plan $g$. 
It is easy to see that, for each couple of plans $f$, $g$, the distance obtained from the $1$-distance is the largest of the four, while the distance obtained from the max-min distance is the smallest.
\noindent
After having defined these distances, we can introduce an interesting variant of the optimization MAPF problem (\ref{op}), namely, the \textit{optimization MAPF problem constrained to a given plan}. This problem is faced when we have a sub-optimal solution $f_0$ of a MAPF instance, and we want to find another solution of the same problem which is not too far from $f_0$ and has better quality, i.e., shorter length. Given $\A^s$ and $\A^t$, initial and final valid configurations on a digraph $G$,  given $f_0 \in \E^*_{\A^s}$ such that $\A^t = \rho(\A^s,f_0) $ (i.e., $f_0$ is a feasible solution of the MAPF instance), and given $r \in {\mathbb{N}}$ and a distance $d$ between plans, the optimization MAPF problem with \textit{Makespan}  constrained to $f_0$ is defined as
 
 \begin{equation} 
 	\label{opp}
 	\centering 
 	\begin{array}{ll}       	
 		
 		\min{|f|} & \\
 		
 		& \\
 		
 		\quad 	\text{s.t.} \;\; \A^t = \rho(\A^s,f) &  \\
 		
 		& \\
 		
 		\quad   \quad  \quad f \in \E^*_{\A^s}, \;  d(f,f_0) \leq r. & \\
 		
 	\end{array}
 \end{equation}

\subsection{Domain reduction of Problems~(\ref{op}) and (\ref{opp})}
\label{eqrel}
In Problems~(\ref{op}) and (\ref{opp}), variable $f$ belongs to the set of well-defined plans $\E^*_{\A^s}$. In order to reduce the cardinality of the feasible set of these two problems, we leverage some invariance properties. Namely, we define two equivalence relations on the set of plans  $\E^*_{\A^s}$ such that the objective function of Problems~(\ref{op}) and (\ref{opp}) has the same value for all plans on the same equivalence class. Further, a plan is feasible if and only if all plans of the same equivalence class are feasible. In this way, we can convert Problems~(\ref{op}) and (\ref{opp}) into equivalent problems that have the set of equivalence classes as the optimization domain. Note that the set of equivalence classes corresponds to the states set that we will use in the dynamic programming solution algorithm.
We will consider the following two equivalence relations on $\E^*_{\A}$. \\

	\begin{defn}
	Let $f_0 \in \E^*_{\A}$ be a reference plan. Given $f,g \in \E^*_{\A}$, then
		\begin{enumerate}
		
		\item $f \sim_1 g$ if and only if 
		
		\begin{enumerate}
			\item $|f| = |g|$;
			
			\item $ \rho(\A,f)=\rho(\A,g)$. 
			
\end{enumerate}

\item $f \sim_2 g$ if and only if 
	
	\begin{enumerate}

 \item $f \sim_1 g$;
		
		\item $d(f,f_0)= d(g,f_0)$. \\

	\end{enumerate}
\end{enumerate}
\end{defn}

We denote by  $\tilde{\E}^i_{\A}$ the set of all equivalence classes of $\sim_i$ on $\E^*_{\A}$. Let $\hat{f} \in \tilde{\E}^i_{\A}$ and   $f \in \E^*_{\A}$ be a representative of the equivalence class $\hat{f}$. We define:
\begin{itemize}
	\item the length of $\hat{f}$, $|\hat{f}| := |f|$;           
	                                        
	\item a new transition function,   
		\[\rho^*: \C \times \tilde{\E}^i_{\A} \rightarrow \C, \]
	\[ \rho^*(\A,\hat{f}):=\rho(\A,f) . \]

	\item the distance from $\hat{f}_0$ (the equivalence class to which $f_0$ belongs), $d(\hat{f},f_0):=d( f,f_0)$ (only if $i=2$).\\
\end{itemize}
            
Note that $|\hat{f}|$ is well-defined, since, by definition, all elements of equivalence class $\hat{f}$ have the same length. Similarly, $\rho^*$ and $d$ are well-defined since, for all elements $f_1,f_2$ of equivalence class $\hat f$, $\rho(\A,f_1)=\rho(\A,f_2)$ and, for $\sim_2$, $d(f_1,f_0)=d(f_2,f_0)$.

Let $\alpha_1: \tilde{\E}^1_{\A}\, \rightarrow {\mathbb{N}} \times \C$ be such that

\begin{equation}\label{a1}
\alpha_1(\hat{f})=(|\hat{f}|,\rho^*(\A,\hat{f}) ).
\end{equation}
This function is well-defined because if $f_1 \sim_1 f_2$ then $|f_1|=|f_2|$ and $\rho(\A,f_1)=\rho(\A,f_2)$. Moreover, $\alpha_1$ is injective because, if $\hat{f}_1$ and $\hat{f}_2$ are such that $\alpha_1(\hat{f}_1) = \alpha_1(\hat{f}_2)$, then $|\hat{f}_1|=|\hat{f}_2|$ and $\rho^*(\A,\hat{f}_1)=\rho^*(\A,\hat{f}_2)$, and, therefore, $\hat{f}_1 = \hat{f}_2$.     
                                                                                                                            
Let $\alpha_2: \tilde{\E}^2_{\A}\, \rightarrow {\mathbb{N}} \times \C \times {\mathbb{N}} $ be defined as follows: 
\begin{equation}\label{a2}
	\alpha_2(\hat{f})=(|\hat{f}|,\rho^*(\A,\hat{f}),d(\hat{f},f_0) ).
\end{equation}

This function is well defined because if $f_1 \sim_2 f_2$, then $|f_1|=|f_2|$, $\rho(\A,f_1)=\rho(\A,f_2)$ and  $d(f_1,f_0)= d(f_2,f_0)$. Moreover, $\alpha_2$ is injective because, if $\hat{f}_1$ and $\hat{f}_2$ are such that $\alpha_1(\hat{f}_1) = \alpha_2(\hat{f}_2)$, then $|\hat{f}_1|=|\hat{f}_2|$, $\rho^*(\A,\hat{f}_1)=\rho^*(\A,\hat{f}_2)$, $d(\hat{f}_1,f_0)= d(\hat{f}_2,f_0)$,  and, therefore, $\hat{f}_1 = \hat{f}_2$.

Since $f \sim_i g$, $i=1,2$, implies that $|f|=|g|$, $\rho(A^s,f)= \rho(A^s,g)$ and $d(f,f_0) \leq r \leftrightarrow d(g,f_0) \leq r$, it turns out that  problems (\ref{op}) and (\ref{opp}) are invariant under the equivalence relations $\sim_i$. Therefore, given $\hat{f}_0\in \tilde{\E}^i_{\A^s}$,  problem (\ref{opp}) (similar for problem (\ref{op})) can be defined as follows over the set of equivalence classes:
\begin{equation} 
 	\label{op3}
 	\centering 
 	\begin{array}{ll}       	
 		
 		\min{|\hat{f}|} & \\
 		
 		& \\
 		
 		\quad 	\text{s.t.} \;\; \A^t = \rho^*(\A^s,\hat{f}) &  \\
 		
 		& \\
 		
 		\quad   \quad  \quad \hat{f} \in \tilde{\E}^i_{\A^s}, \;  d(\hat{f},f_0) \leq r. & \\
 		
 	\end{array}
 \end{equation}

\subsection{Neighborhoods}
Given the distances defined in Section \ref{sec:distances} and the definition of the equivalence classes in Section \ref{eqrel}, we can define the neighborhood of a vertex, of a configuration, and of an equivalence class. Moreover, we are able to upper estimate the cardinality of such neighborhoods.
Such estimate is needed to evaluate the time needed to explore the neighborhoods, an operation that is central in the approach proposed in this paper.\\
\noindent
First, we define a ball $\B_r(v)$ centered at $v \in V$ of radius $r \in {\mathbb{N}}$:
\[\B_r(v) := \{u \in V: \; d(u,v)\leq r\}.\]
We denote by $\bar{\B}_r(v)$ the border of the ball, obtained by replacing $\leq$ with $=$ in the definition.
Let $\phi =outdeg(G)$ be the maximum out-degree of digraph $G=(V,E)$.
The following proposition provides an upper bound on the cardinality of $\bar{\B}_r(v)$.
\begin{prop}
It holds that
\begin{equation} \label{eq2}
		|\bar{\B}_r(v)| \leq  \phi^r.
\end{equation}
\end{prop}	
\textit{Proof.} Let  $n_h$ be the number of nodes at distance $h$ from $v$.  Note that $n_1 \leq \phi$, and $\forall h \geq 2$ $n_h \leq n_{h-1} (\phi -1)$. By induction, $n_h\leq \phi (\phi -1)^{h-1}$ $\forall h \geq 1$. Therefore, an upper-bound for the number of nodes on the border of the ball is 
\[|\bar{\B}_r(v)| \leq \phi (\phi -1)^{r-1} \leq \phi^r.\]
\noindent
Next, we define a ball centered at $\A \in \C$ of radius $r$:
\[\B_r(\A) := \{\A^* \in \C: \; d(\A^*,\A)\leq r\},\]
denoting by $\bar{\B}_r(\A)$ its border.
An upper bound for the cardinality of the ball is given by the following proposition.
\begin{prop}
It holds that
\begin{equation}\label{eq3}
		|\B_r(\A)| \leq 1 + \frac{(r +k-1)!}{(r-1)! (k-1)!}  \phi^r.
	\end{equation}
\end{prop}
\textit{Proof.} First of all, we find an upper bound for the number of configurations at distance $h$ from $\A$, i.e., an upper-bound for the cardinality of the border of the ball centered in $\A$ of radius $h$:

\[ \bar{\B}_h(\A) = \left\{\A^* \in \C: \;\sum_{i=1}^{|P|} d(\A(p_i),\A^*(p_i)) = h\right\}. \]

Let $(h_1, \cdots,h_{|P|})$ be a $|P|$-decomposition of $h$ (i.e., $\sum_{i=1}^{|P|} h_i = h $ ). An upper bound for the number of configurations for which $ d(\A(p_i),\A^*(p_i)) = h_i$ is
\[ \prod_{i=1}^{|P|}|\bar{\B}_{h_i}(\A(p_i))| \leq \prod_{i=1}^{|P|} \phi^{h_i}  = \phi^h. \]
By standard combinatorial arguments, the number of $|P|$-decompositions of $h$ is 
\[   \frac{(h + (|P|-1))!}{h! (|P|-1)!} .\]

Therefore, the cardinality of the border of the ball of radius $h$ can be bounded from above by:
\[|\bar{\B}_{h}(\A)| =  \frac{(h + (|P|-1))!}{h! (|P|-1)!} \prod_{i=1}^{|P|}|\bar{\B}_{h_i}(\A(p_i))| \leq\]
\[\leq  \frac{(h + (|P|-1))!}{h! (|P|-1)!} \phi^h,\]
and the total number of configurations in $\B_{r}(\A)$ can be overestimated as follows: 
\[ |\B_{r}(\A)| = 1 +  \sum_{h=1}^{r}|\bar{\B}_{h}(\A)| \leq \]
\[  \leq 1 +\sum_{h=1}^{r}  \frac{(h +|P|-1)!}{h! (|P|-1)!}  \phi^h \leq 1 + r \,\frac{(r +|P|-1)!}{r! (|P|-1)!}  \phi^r.\]
\noindent
Finally, given a radius $r \in {\mathbb{N}}$, we define the following \textbf{neighborhood} of $f_0\in {\E}^*_{\A}$:
\[\N_r(f_0) := \{\hat{g} \in  \tilde{\E}^i_{\A}: \; |\hat{g}| \leq |f_0|, \,d(\hat{g},f_0)\leq r\}.\]\label{ball}
Here we consider the distance based on the max-min distance. 
The following proposition provides an upper bound on the cardinality of $\N_r(f_0)$. Note that, since the max-min distance is the smallest among the considered distances, the upper bound provided in the proposition is also valid for the distances
based on the other distances previously discussed. 
\begin{prop}
	\label{neig3}
It holds that:
	\begin{equation}\label{eq4}
		|\N_r(f_0)| \leq   \,|f_0|^2 \,\left(  1 + \frac{(r +k-1)!}{(r-1)! (k-1)!}  \phi^r\right),
	\end{equation} 
where $k=|P|$ and $\phi=outdeg(G)$.
\end{prop}
\textit{Proof.}
Let $\hat{f}\in  \tilde{\E}^i_{\A}$, $ f_0\in  {\E}^*_{\A}$ with
$\hat{f}\in \B_r(f_0)$.
 Let $f$ be a representative of $\hat{f}$ and let $I = \{1, \cdots, |\hat{f}|\}  $ and $J = \{1, \cdots, |f_0|\}  $. Then:
\[\alpha_1(\N_r(f_0)) \subset  \left(I \times  \bigcup_{j \in J} \B_r(\psi_{f_0}(j)) \right).\]
		        Indeed, for each representative $f$ of $\hat{f}$, it holds that $|f| \leq |f_0|$ and, moreover, $d(\hat{f},f_0)\leq r$ implies that 
$$
\max_{1 \leq i \leq |f|} \min_{1 \leq j \leq |f_0|} d(\psi_{f}(i), \psi_{f_0}(j)) \leq r,
$$ 
and, in particular, for $i=|f|$ there exists $j$ such that
\[d(\rho(\A, f), \psi_{f_0}(j)) \leq r, \]
which means that $\rho(\A, f) \in \B_r(\psi_{f_0}(j))$.  Therefore, recalling that $\alpha$ is injective,

	        \[|\N_r(f_0)| \leq \big|I \times  \bigcup_{j \in J} \B_r(\psi_{f_0}(j)) \big| .\]
Then, also in view of  (\ref{eq3}), we have that:
\[ |\N_r(f_0) | \leq  |f_0| \left( \sum_{j=1}^{|f_0|}|\B_{r}(\psi_{f_0}(j))| \right) \leq \]

    \[ \,\leq |f_0|^2 \,\left(  1 + \frac{(r +k-1)!}{(r-1)! (k-1)!}  \phi^r\right).\]		

\begin{prop}
	\label{neig4}
	The neighborhood of $f_0$ of radius $r$ has a polynomial cardinality with respect to the number of nodes. In particular, $\exists \, C=C(r)\in {\mathbb{R}}$ such that	
	\[|\N_r(f_0) | \leq |f_0|^2 \left( 1+ C (r+k)^r\phi^r \right), \]	
	where $k=|P|$ and $\phi=outdeg(G)$.
	
	
\end{prop}
\textit{Proof.} Note that
\[\frac{(r +k-1)!}{(r-1)! (k-1)!} \leq \frac{(r + k)^r}{(r-1)!}.\]
Defining $C:=\frac{1}{(r-1)!}$ and reminding the result of Proposition~\ref{neig3} we have the thesis.

\section{Iterative Neighborhood Search}
In this section we describe the proposed iterative approach to detect a sub-optimal solution of Problem~(\ref{op}), or the equivalent counterpart of this problem defined over the set of equivalence classes. The returned solution is locally optimal with respect to the neighborhood of a reference solution. At each iteration, we solve an instance of Problem~(\ref{op3}).
More in detail, the algorithm takes as input a feasible solution $f_0$, that may be of poor quality. For instance, we can obtain $f_0$ from a rule-based algorithm. We aim at improving $f_0$, obtaining a shorter solution.
To this end, we solve Problem~(\ref{op3}) with a dynamic programming algorithm. Namely, we search in neighborhood $\N_r(f_0)$ for plans shorter than $f_0$ through algorithm \emph{Dynprog}, that we will describe below.
If we cannot obtain a solution shorter than $f_0$ (that is, $f_0$ is locally optimal) we stop the algorithm. Otherwise, if we obtain an improved solution $f^*$, we redefine the reference solution as $f_0=f^*$ and solve again Problem~(\ref{op3}). We iterate this procedure until we cannot shorten the current solution any further.

This algorithm can be classified as a \textbf{Neighborhood Search} algorithm (see~\cite{pisinger2019}).

To define the neighborhood $N_r(f_0)$, we can use any distance function among those presented in Section~\ref{sec:distances}. In our numerical experiments, we used the sum-min distance. 
Algorithm~\ref{neigh} presents the steps of the procedure. 

\begin{algorithm}
\caption{\label{neigh}Neighborhood Search}
\label{alg1}

\begin{algorithmic}
\State Input: $f_0, r, \A^s, \A^t$
\State Output: $f^*$
\State $f^* \gets f_0$

\Do
\State $f_0 \gets f^*$
    \State $f^* \gets DynProg(f_0, r, \A^s, \A^t)$
\doWhile{$|f^*| < |f_0|$} 

\State \textbf{return} $f^*$
 
\end{algorithmic}
\end{algorithm}

In Algorithm~\ref{neigh}, $r$ is the local search radius, $\A^s$ is the initial configuration, and $A^t$ is the final one.

\section{Dynamic Programming Algorithm}


To search for the optimal solution of Problem (\ref{op3}), we employ a \textit{Dynamic Programming} (DP) algorithm. In generic DP problems we are given a state space $S$ where $A \subset S$ are the target states, an expansion function $g: S \rightarrow P(S)$, where $P$ is the power set of $S$, and an objective function $c : S \rightarrow R$. Starting from an initial state $s_0 \in S$, we iteratively expand states with function $g$ to explore the state space and compute $s_t = g^n(s_0) \in A$ with the minimal objective function.
In our case, the states represent the equivalence classes of relation $\sim_2$, defined in Section~\ref{eqrel}. Namely, we use injection $\alpha_2: \tilde{\E}^2_{\A}\, \rightarrow {\mathbb{N}} \times \C \times {\mathbb{N}}$ to associate to each equivalence class $\hat f$ a triple $(\beta,\gamma,\sigma)=\alpha_2(\hat f)$, where $\beta$ is the length of $\hat f$, $\gamma$ the configuration obtained by applying a plan representative of $\hat f$ to the initial state $\A^s$, and $\sigma$ is the distance of a representative of $\hat f$ from reference plan $f_0$.
Namely, the state space is  
	\[\Ss := \alpha_2(\tilde{\E}^2_{\A} ) \subset {\mathbb{N}} \times \C \times {\mathbb{N}},\]
	where $\alpha_2$ is defined in (\ref{a2}). Since $\alpha_2$ si injective, $\Ss $ and $\tilde{\E}^2_{\A}$ are in one-to-one correspondence.
Each state $s= (\beta, \gamma, \sigma) \in \Ss$, represents the equivalence class:
\[   \alpha_2^{-1}(s) = \{ f \in \E_{\A} \, : \, \beta = |f|, \,\gamma = \rho(\A^s, f), \]
\[\,\sigma = d(f,f_0) \}.\]  
The initial state is $s_0=\alpha_2(\epsilon)=(0, \A^s, 0)$.
We use a priority queue $Q$ to store the states that have not been visited yet. At the beginning, $Q=\{s_0\}$.
We define a partial ordering on $\Ss$ based on length. Namely, if $s_1=(\beta_1,\gamma_1,\sigma_1),s_2=(\beta_2,\gamma_2,\sigma_2)$, $s_1 < s_2$ if $\beta_1 < \beta_2$. We order the elements of $Q$ according to this partial ordering.

A state $s_1=(\beta_1,\gamma_1,\sigma_1)$ \emph{dominates} $s_2=(\beta_2,\gamma_2,\sigma_2)$ if
\begin{itemize}
\item $\beta_1 \leq \beta_2$,
\item $\gamma_1=\gamma_2$,
\item $\sigma_1 \leq \sigma_2$.
\end{itemize}

In other words, $s_1$ dominates $s_2$ if the plans $f_1$, $f_2$, corresponding to $s_1$ and $s_2$, satisfy the following properties. Plan $f_1$ is not longer than $f_2$, $f_1$ and $f_2$ lead to the same final configuration, and the distance of $f_1$ from the reference solution $f_0$ is not larger than the one of $f_2$. If $s_1$ dominates $s_2$, we can discard $s_2$. In general, we remove from $Q$ all dominated states.
\noindent
We also define the following {\em transition function}, which allows to (possibly) add new states to the priority queue: 
\[\tilde{\rho}:\Ss \times \E \rightarrow \Ss  \]
\[ \tilde{\rho}((\beta,\gamma,\sigma), e) :=  (\beta + 1, \rho(\gamma, e), \sigma +   \min_{k \in {\mathbb{N}}} d(\rho(\gamma, e), \psi_{f_0}(k))). \]
Applying this function on a state $s=(\beta, \gamma, \sigma)$:
\begin{itemize}
\item adds $1$ to the length of the class $\alpha_2^{-1}(s)$;
\item updates the final configuration of the equivalence class, applying the function $\rho(\gamma, e)$ to the final configuration of $\alpha_2^{-1}(s)$ with $e$ being the chosen set of edges;
\item updates $\sigma$, adding the computed minimum distance between the updated final configuration $\rho(\gamma, e)$ and the reference plan $f_0$.

\end{itemize}

We denote by
\[\Sigma:= \left\{\tilde{\rho}((\beta,\gamma,\sigma),e)\ :\  e \in \E\  \mbox{and}\ \rho(\gamma, e) \in \B_r(\gamma)\right\}\subset \Ss,\]
the set of new states which can be generated through the transition function $\tilde{\rho}$ applied to the current state $(\beta,\gamma,\sigma)$ and all possible actions in $\E$ leading to configurations in $\B_r(\gamma)$.  Moreover, we denote with 
\[  \Gamma := \alpha_2(\N_r(f_0)) \subset \Ss,\]
the set of states that can be visited during a neighborhood search.
\subsection{Algorithm}
The Dynamic Programming algorithm is described in Algorithm \ref{dp_alg}.
The priority queue $Q$ maintained inside the algorithm is a set of states, ordered by the length $\beta$ of their representatives. Function ${\tt insert}(Q,x)$ inserts a state $x$ maintaining the partial order of $Q$.
Function ${\tt remove}(Q,x)$ removes $x$ from $Q$.  The head of the queue, that is the state with minimal $\beta$, is denoted by $Q[0]$. The algorithm explores the state space starting from the initial state $s_0$. At each iteration, the state with minimum $\beta$ is extracted from the queue. If the state extracted is the target state (that is, if $\gamma=\A^t$) the algorithm stops and we 
return a representative of the optimal solution of Problem  (\ref{op3}) (for a given equivalence class $\hat{f}$, the function ${\tt repr}(\hat{f})$ returns a representative of the class).
Otherwise, the algorithm employs function ${\tt expand}(s,f_0,r)$, based on
the transition function previously defined, to find new states. If a new state is not dominated, then it is added to the queue $Q$. Moreover, all states in $Q$ dominated by the newly added state are removed from $Q$.
\begin{algorithm}
\caption{Dynamic Programming with Dominance}
\label{dp_alg}
\begin{algorithmic}

\State Input: $f_0, r, \A^s, \A^t$
\State Output: $f$
\BState $s_0 \gets (0, \A^s, 0)$
\State ${\tt insert}(Q, s_0)$
\While{$Q \neq \emptyset$}:

	\State $s=(\beta,\gamma,\sigma) \gets Q[0]$
	\If{$\gamma=\A_t$}
		\State $f \gets {\tt repr}(\alpha_2^{-1}(s))$   
\State $Q \gets \emptyset$
 \Else
	\State $\Sigma \gets {\tt expand}(s, f_0, r)$
	\For{$s_k \in \Sigma$}
		\If{$s_k$ is not dominated in $Q$}
			\State ${\tt insert}(Q, s_k)$ 
			\For{$s_i \in Q$}
				\If{$s_k$ dominates $s_i$}
					\State ${\tt remove}(Q, s_i)$
				\EndIf	
			\EndFor
		\EndIf
	\EndFor
\EndIf
\EndWhile
\State
\State \textbf{return} $f$
\end{algorithmic}
\end{algorithm}
Note that for more complicated solutions, the algorithm is far more effective. 
\begin{theorem}
  \label{thm_complx}
Algorithm \ref{alg1} and \ref{dp_alg} have polynomial time complexity with respect to the number of nodes of the graph.
\end{theorem}
\textit{Proof.} In Algorithm \ref{dp_alg}, the time complexity is  $O(|Q|^2 \cdot |\Sigma|)$.
Sets $Q$ and $ \Sigma$ vary with each iteration, but always remain subsets of $ \Gamma$. So at each iteration the following upper bound for their cardinality always holds: 
\[|Q|  \leq |\Gamma| , |\Sigma| \leq |\Gamma|  .\]
Reminding that $\alpha_2$ is injective and using the result of Proposition \ref{neig4},
\[ |\Gamma|= |\N_r(f_0)|\leq |f_0|^2 \left( 1+ C (1+k)^r\phi^r \right), \]	
where $k=|P|$ and $\phi=outdeg(G)$. Therefore, the time complexity of Algorithm \ref{dp_alg} is
\[O(|f_0|^6 \left( 1+ C (r+k)^r\phi^r \right)^3).\]
Algorithm \ref{alg1} recalls at most $|f_0|$ times Algorithm \ref{dp_alg}, and so it has time complexity $O(|f_0|^7 \left( 1+ C (r+k)^r\phi^r \right)^3)$.

\section{Experimental results}
We performed two distinct sets of experiments on different graphs and with initial solutions generated in two distinct ways.
The algorithm has been coded with the \textit{C++} programming language, and has been run on a \textit{11th Gen Intel(R) Core(TM) i7-1165G7 @ 2.80GHz} processor with a 16 GB RAM.

\subsection{Random Graphs with sequentially generated initial solution}
\begin{figure}[h]
\centering
\includegraphics[width=0.32\textwidth]{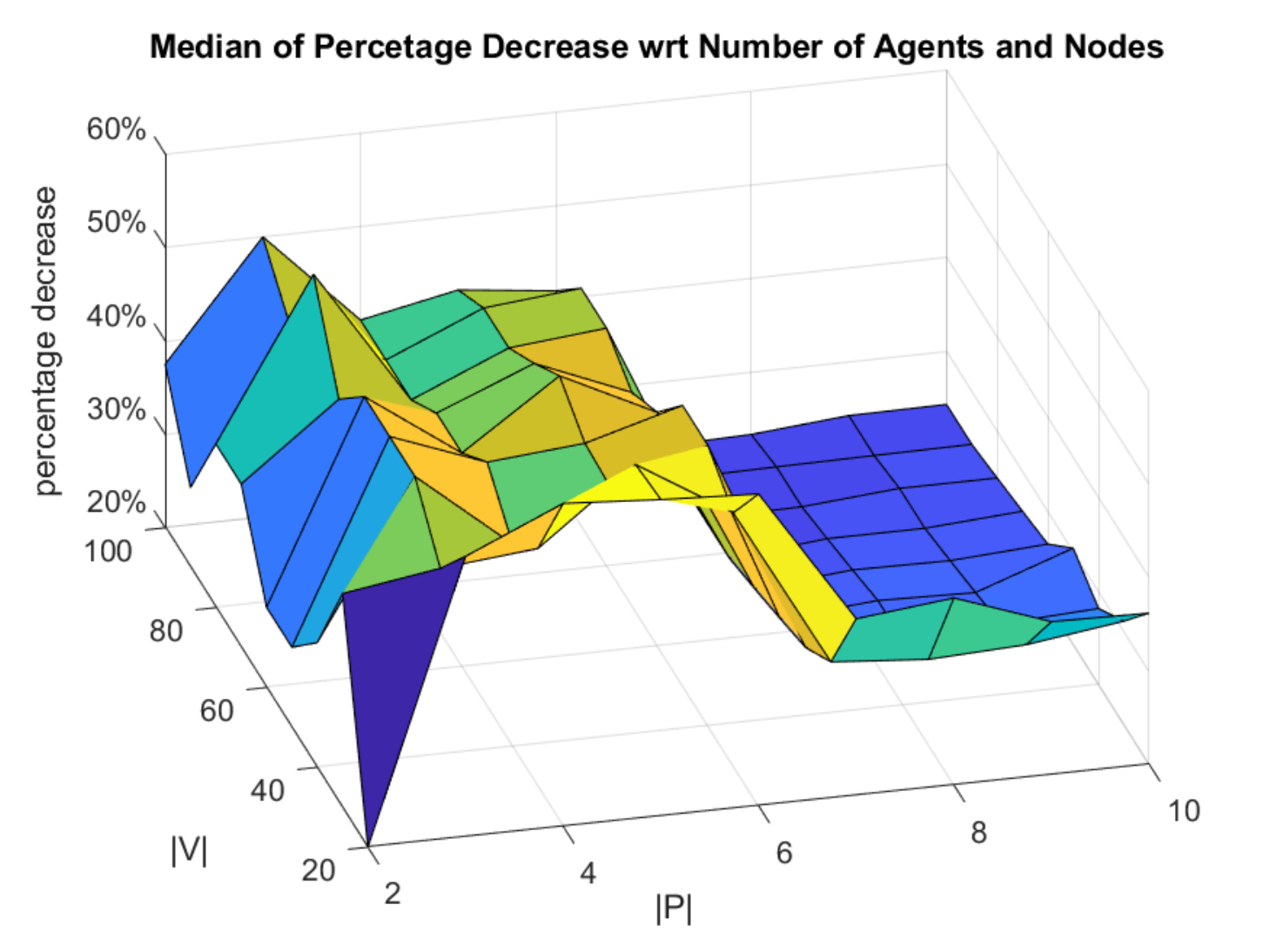}
\caption{Average Percentage Decrease per $|P|$ and $|V|$.}
\label{fig:pd}
\end{figure}

\begin{figure}[h]
\centering
\includegraphics[width=0.32\textwidth]{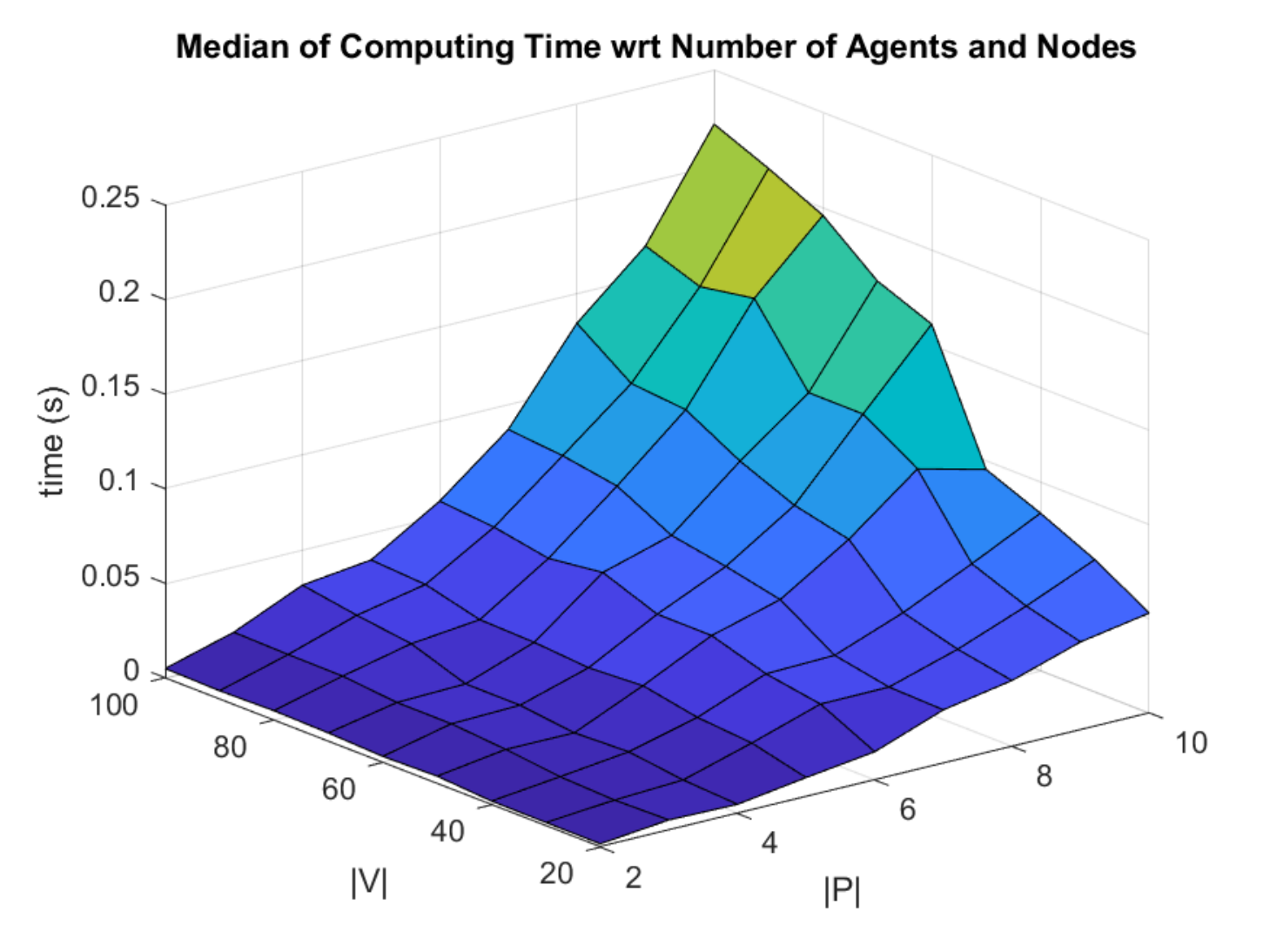}
\caption{Average Running Time per $|P|$ and $|V|$.}
\label{fig:act}
\end{figure}

In the first set of experiments we generated random directed graphs with a number of nodes $|V|$ ranging from 20 to 100 by 10, and a number of edges $|E|$ equal to $4|V|$. The graphs are generated by creating $|E|$ random ordered pair of nodes and using them to build a directed graph. Only strongly connected graphs are selected. The number of agents $|P|$ ranges from 2 to 10, while $\A^s$ and $\A^t$ are randomly generated. In these experiments, to generate the initial solutions, each agent is brought to its target node one at a time, following the shortest path, in terms of crossed edges, from its source to its target in the graph obtained by removing the nodes currently occupied by all other agents (either their target or their source, depending on whether they have been already moved or not). 
Note that such procedure is not complete, i.e., it does not guarantee to find a feasible solution when it exists or to establish that no such feasible solution exists.
We generated 100 random graphs, for which the described procedure was able to return a feasible solution, for every combination of number of nodes and number of agents. 
After some tuning, we set the radius $r$ of the neighborhood equal to 5, which turned out to be a good compromise between the quality of the solutions found in the neighborhood and the time needed to explore the neighborhood (note that
such time increases exponentially with $r$).
Given the initial solution $f_0$ and the final one $f^*$ returned by the proposed approach, the percentage decrease of the final solution wrt the initial solution is equal to $100\frac{f_0-f^*}{f_0}\%$. 
In Figures \ref{fig:pd} and \ref{fig:act} we report the median of the average percentage decrease and of the running time (in seconds), respectively, for every combination of $|P|$ and $|V|$. 
It is worthwhile to remark that the percentage decrease tends to be lower as the number of agents increases. A tentative explanation is that for cases with a greater number of agents, when the sequential procedure to generate an initial solution is able to return a feasible solution, such solution is already a good one which cannot be largely improved. This phenomenon is not observed in the second set of experiments, where a different procedure to generate an initial feasible solution is employed.

Note that the average percentage decrease is lower when the number of agents becomes higher. This can be explained by the nature of the algorithm employed in finding the initial solution. Increasing the number of agents without increasing the number of nodes gives us way less feasible instances and the solutions found are more difficult to improve.


\subsection{Strongly connected with multiple biconnected components graphs and rule-based generated initial solution}

\begin{figure}[h]
\centering
\includegraphics[width=0.32\textwidth]{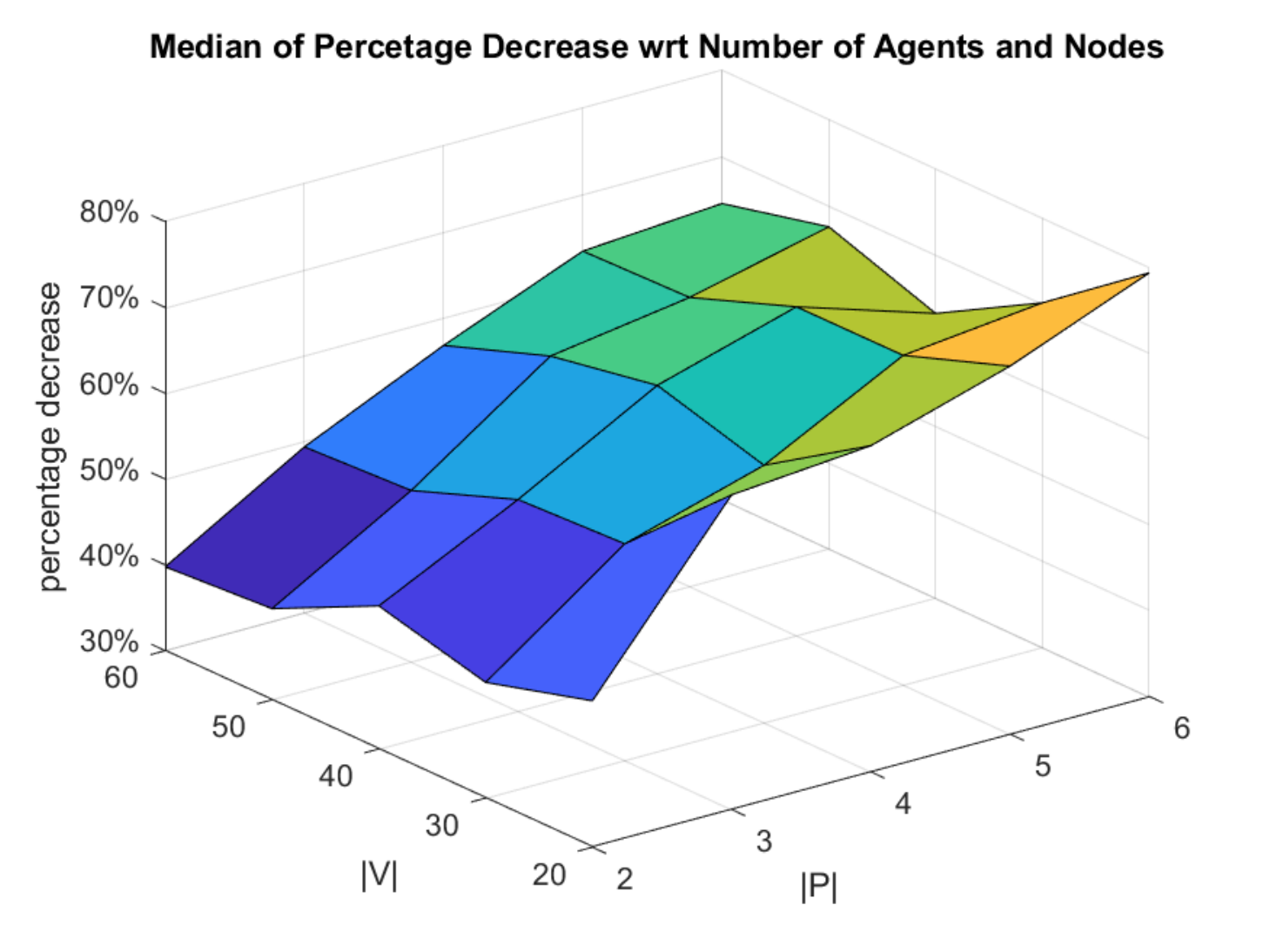}
\caption{Average Percentage Decrease per $|P|$ and $|V|$.}
\label{fig:disc_pd}
\end{figure}

\begin{figure}[h]
\centering
\includegraphics[width=0.32\textwidth]{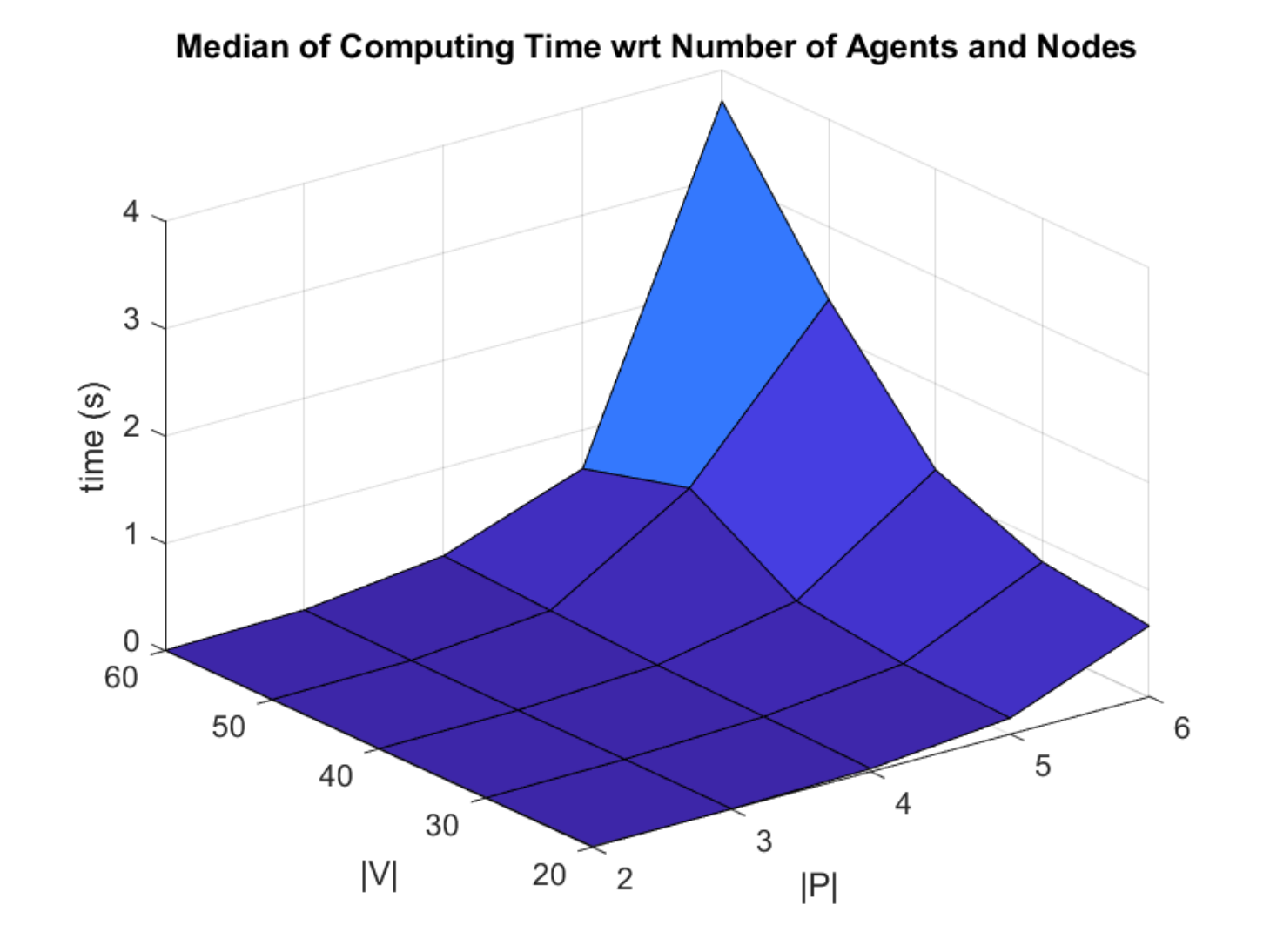}
\caption{Average Running Time per $|P|$ and $|V|$.}
\label{fig:disc_act}
\end{figure}

 In the second set of experiments, we generated graphs that are strongly connected and with multiple biconnected components. The procedure used to generate such graphs can be found in \cite{ardizzoni2022}.
The number of agents $|P|$ ranges from 2 to 6, and the number of nodes vary from 20 to 50 with an increment of 10. The initial and final configuration $\A^s$ and $\A^t$ are randomly generated. 
In this case the initial solutions (if they exist) are generated through the algorithm described in \cite{ardizzoni2022}. Note that this algorithm is complete, i.e., it always returns a feasible solution in case one exists. However, the quality of such initial solution
is usually lower (i.e., the initial plan is usually longer) with respect to the one returned by the algorithm employed in the first set of experiments to generate an initial solution. This might be also the explanation why the percentage decrease in these experiments (see Figure \ref{fig:disc_pd}) appears to be larger wrt the first set of experiments, and also tends to increase with the number of agents (differently from the case of random graphs).
Again after some tuning, we set the radius $r$ of the neighborhood equal to 3.
As for the first set of experiments, Figures \ref{fig:disc_pd} and \ref{fig:disc_act} report the median of the average percentage decrease and of the running time (in seconds), respectively, for every combination of $|P|$ and $|V|$.
Note that the graphs employed in the second set of experiments appear to be more challenging with respect to the random ones. Indeed, computing times are larger and increase rapidly with the number of agents.


\section{Conclusion and future works}
We proposed an iterative local search procedure for MAPF, in order to shorten a known feasible solution. We obtain a solution, that is still sub-optimal, but, in general, of much better quality than the initial one.
The proposed algorithm has polynomial time complexity with respect to the number of agents (see Theorem~\ref{thm_complx}) and has computational times compatible with industrial applications. We can extend the results in various respects, that will be the focus of future research:
\begin{itemize}
\item We can define locality constraints different from the ones considered in Section~\ref{sec:distances}. For instance, we can set a maximum on the number of agents that modify their path with respect to the reference solution. Alternatively, we can set an upper bound on the number of time intervals in which the solution departs from the reference one.
\item We can test the proposed approach on real industrial scenarios.
  \item We can improve the complexity bound presented in Theorem~\ref{thm_complx}, that it currently based on a quite rough bound.
  \end{itemize}

\printbibliography 
\end{document}